\documentstyle[amscd,amssymb,verbatim,12pt]{amsart}

\theoremstyle{plain}
\newtheorem{Thm}[subsection]{Theorem}

\newtheorem{Cor}[subsection]{Corollary}

\newtheorem{Lem}[subsection]{Lemma}
\newtheorem{Prop}[subsection]{Proposition}

\theoremstyle{definition}
\newtheorem{Def}[subsection]{Definition}

\theoremstyle{remark}

\newtheorem{Rem}{Remark}

\renewcommand{\rm}{\normalshape}

\newif\ifShowLabels
\ShowLabelstrue
\newdimen\theight
\def\TeXref#1{%
	\leavevmode\vadjust{\setbox0=\hbox{{\tt
		\quad\quad  {\small \rm #1}}}%
	\theight=\ht0
	\advance\theight by \lineskip
	\kern -\theight \vbox to
	\theight{\rightline{\rlap{\box0}}%
	\vss}%
	}}%

\ShowLabelsfalse

\newcommand{\ssec}[2]{\subsection{#2}\label{SS:#1}%
	\ifShowLabels \TeXref{{SS:#1}} \fi}

\newcommand{\refss}[1]{Section ~\ref{SS:#1}}

\newcommand{\reft}[1]{Theorem ~\ref{T:#1}}

\newcommand{\refe}[1]{\eqref{E:#1}}

\newenvironment{thm}[1]%
	{ \begin{Thm} \label{T:#1}  \ifShowLabels \TeXref{T:#1} \fi }%
	{ \end{Thm} }

\newcommand{\th}[1]{\begin{thm}{#1} }
\renewcommand{\eth}{\end{thm} }

\newenvironment{lemma}[1]%
	{ \begin{Lem} \label{L:#1}  \ifShowLabels \TeXref{L:#1} \fi }%
	{ \end{Lem} }
\newcommand{\lem}[1]{\begin{lemma}{#1}}
\newcommand{\elem}{\end{lemma}}

\newenvironment{propos}[1]%
	{ \begin{Prop} \label{P:#1}  \ifShowLabels \TeXref{P:#1} \fi }%
	{ \end{Prop} }
\newcommand{\prop}[1]{\begin{propos}{#1}}
\newcommand{\eprop}{\end{propos}}

\newenvironment{corol}[1]%
	{ \begin{Cor} \label{C:#1}  \ifShowLabels \TeXref{C:#1} \fi }%
	{ \end{Cor} }
\newcommand{\cor}[1]{\begin{corol}{#1}}
\newcommand{\ecor}{\end{corol}}

\newenvironment{defeni}[1]%
	{ \begin{Def} \label{D:#1}  \ifShowLabels \TeXref{D:#1} \fi }%
	{ \end{Def} }
\newcommand{\defe}[1]{\begin{defeni}{#1}}
\newcommand{\edefe}{\end{defeni}}

\newenvironment{remark}[1]%
	{ \begin{Rem} \label{R:#1}  \ifShowLabels \TeXref{R:#1} \fi }%
	{ \end{Rem} }
\newcommand{\rem}[1]{\begin{remark}{#1}}
\newcommand{\erem}{\end{remark}}

\newcommand{\eq}[1]%
	{ \ifShowLabels \TeXref{E:#1} \fi
	   \begin{equation} \label{E:#1} }
\newcommand{\eeq}{ \end{equation} }

\newcommand{\prf}{ \begin{pf} }
\newcommand{\epr}{ \end{pf} }

\newcommand\alp{\alpha}		
\newcommand\bet{\beta}
		\newcommand\Gam{\Gamma}
		\newcommand\Del{\Delta}

\newcommand\lam{\lambda}

\newcommand\ome{\omega}		\newcommand\Ome{\Omega}

\newcommand\calF{{\cal{F}}}

\newcommand\calL{{\cal{L}}}
\newcommand\calM{{\cal{M}}}
\newcommand\calN{{\cal{N}}}

\newcommand\calP{{\cal{P}}}
\newcommand\calQ{{\cal{Q}}}

\newcommand\RR{\Bbb{R}}

\newcommand\CC{\Bbb{C}}

\newcommand\nek{,\ldots,}
\newcommand\sdp{\times \hskip -0.3em {\raise 0.3ex
\hbox{$\scriptscriptstyle |$}}} 

\newcommand\ind{\operatorname{ind}}

\newcommand\Ker{\operatorname{Ker}}

\newcommand{\F}{\calF}

\newcommand{\w}{\text{\( \ome\)}}

\renewcommand{\t}{\text{\( \tau_{\sqrt{t}}\)}}

\newcommand{\n}{\nabla}

\renewcommand{\b}{\bullet}

\begin{document}


\title[]{The Novikov-Bott inequalities}
\author[]{Maxim Braverman} \author[]{Michael Farber}
\address{School of Mathematical Sciences\\
Tel-Aviv University\\
Ramat-Aviv 69978, Israel}
\email{maxim@@math.tau.ac.il, \,farber@@math.tau.ac.il}
\thanks{The research was supported by grant No. 449/94-1 from the
Israel Academy of Sciences and Humanities.}

{\ }
\vspace{0.7cm}

\maketitle

\vspace{-0.6cm}
\begin{abstract}
We generalize the Novikov inequalities for 1-forms in two different
directions: first, we allow non-isolated critical points (assuming
that they are non-degenerate in the sense of R.Bott), and, secondly,
we strengthen the inequalities by means of twisting by an arbitrary
flat bundle.  We also obtain an $L^2$ version of these inequalities
with finite von Neumann algebras.

The proof of the main theorem uses Bismut's modification of the Witten
deformation of the de Rham complex; it is based on an explicit
estimate on the lower part of the spectrum of the corresponding
Laplacian.
\end{abstract}
\vspace{-0.6cm}

\ssec{novnum}{The generalized Novikov numbers}Let $M$ be a closed manifold and
let $\F$ be a complex flat vector bundle over $M$. We will denote by
$\n: \Ome^\b(M,\F)\to \Ome^{\b+1}(M,\F)$ the covariant derivative on $\F$.
Given a closed 1-form $\omega\in \Omega^1(M)$ on $M$ with real values,
it determines a family of connections on $\F$ (the {\em Novikov deformation})
parameterized by the real numbers $t\in \RR$
\eq{nab-t}
	\nabla_t:\Omega^i(M,\F)\to \Omega^{i+1}(M,\F),\qquad
		\n_t:\, \theta\mapsto \n\theta+t\ome\wedge\theta.
\end{equation}

All the connections $\nabla_t$ are flat, i.e. $\nabla_t^2=0$.  Denote
by $\F_t$ the flat vector bundle defined by the connection
\refe{nab-t}. Note that changing $\ome$ by a cohomologious 1-form
determines a gauge equivalent connection $\nabla_t$ and so the
cohomology $H^\b(M,\F_t)$ depends only on the cohomology class $\xi=
[\omega]\in H^1(M,\RR)$ of $\omega$.  One can show that there exists a
{\em finite} subset $S\subset\RR$ (the set {\em jump points}) such that
$\dim H^i(M,\F_t)$ \/ is constant for $t\notin S$ and it jumps up for
$t\in S$. The dimension of $\dim H^i(M,\F_t)$ \/ for $t\notin S$ \/ is
called the $i$-th {\em (generalized) Novikov number}
$\beta_i(\xi,\F)$.

\ssec{assum}{Assumptions on the 1-form} Let $C$ denote the set of
critical points of \w\ (i.e. the set of points of $M$, where \w\
vanishes).  We assume that \w\ is {\em non-degenerate in the sense of
Bott}, i.e.  $C$ is a submanifold of $M$ and that the Hessian of \w\
is a non-degenerate quadratic form on the normal bundle $\nu(C)$ to
$C$ in $M$.  Here by the Hessian of \w\ we understand the Hessian of
the unique function $f$ defined in a tubular neighborhood of $C$ and
such that $df=\w$ and $f_{|_C}=0$.

\ssec{main}{The main result} Let $Z$ be a connected component of the
critical point set $C$ and let $\nu(Z)$ denote the normal bundle to
$Z$ in $M$. Since the Hessian of \w\ is non-degenerate, the
bundle $\nu(Z)$ splits into the Whitney sum of two subbundles
$\nu(Z) = \nu^+(Z) \oplus \nu^-(Z)$,
such that the Hessian is strictly positive on $\nu^+(Z)$ and strictly
negative on $\nu^-(Z)$. The
dimension of the bundle $\nu^-(Z)$  is called the {\em index}
of $Z$ (as a critical submanifold of \w) and is denoted by $\ind(Z)$.
Let $o(Z)$ denote
the {\it orientation bundle of $\nu^-(Z)$, considered as a flat
line bundle}. Consider the {\it twisted Poincar\'e polynomial} of $Z$
\eq{Poin}
	\calP_{Z,\F}(\lambda)\ =\
		\sum \lambda^i\dim_{\CC} H^i(Z,\F_{|_Z}\otimes o(Z))
\end{equation}
(here $H^i(Z,\F_{|_Z}\otimes o(Z))$ denote the cohomology of $Z$ with
coefficients in the flat vector bundle $\F_{|_Z}\otimes o(Z)$)
and define using it the following {\it Morse counting polynomial}
\eq{M-pol}
	\calM_{\omega,\F}(\lambda)\ =
		\ \sum_{Z} \lambda^{\ind(Z)} \calP_{Z,\F}(\lambda),
\end{equation}
where the sum is taken over all connected components $Z$ of $C$.

With one-dimensional cohomology class
$\xi=[\omega]\in \-H^1(M,\RR)$ and the flat vector bundle $\F$, one can
associate the {\em Novikov  polynomial}
\eq{N-pol}
	\calN_{\xi,\F}(\lambda)\ =\ \sum_{i=0}^n \lambda^i
		\beta_i(\xi,\F),  \qquad n=\dim M.
\end{equation}

\th{main} There exists a polynomial
  $\calQ(\lam)$ with non-negative integer coefficients, such that
  \eq{main}
	\calM_{\omega, \F}(\lambda)\ -\ \calN_{\xi,\F}(\lambda)\ =\
			(1+\lambda)\calQ(\lambda).
  \end{equation}
\eth

The main novelty in this theorem is that it is applicable to the case
of 1-forms with non-isolated singular points. Thus, we obtain, in
particular, a new proof of the degenerate Morse inequalities of
R.Bott.  Moreover, \reft{main} yields a generalization of the
Morse-Bott inequalities to the case of an arbitrary flat vector bundle
$\F$; this generally produces stronger inequalities  as shown in
\refss{example}.


\begin{Cor}[Euler-Poincar\'e theorem] \label{C:Eu-Poin}
  Under the conditions of \reft{main}, the Euler
  characteristic of $M$ can be computed as
  $
	\chi(M)\ = \sum_{Z}  (-1)^{\ind(Z)}\chi(Z),
  $
  where the sum is taken over all connected components $Z\subset C$.
\end{Cor}

\ssec{isol}{The case of isolated critical points}
Consider now the special case when all critical points of $\omega$
are isolated.  \reft{main} gives  the inequalities
\eq{novin''}
	\sum_{i=0}^{p}(-1)^im_{p-i}(\omega)\
		\ge\  d^{-1}\cdot \sum_{i=0}^{p}(-1)^i\beta_{p-i}(\xi, \F),
				\qquad p=0,1,2,\dots,
\end{equation}
where $d=\dim \F$ and $m_p(\omega)$ denotes the number of critical points
of $\omega$ of index $p$.
The last inequalities coincide with the Novikov inequalities
\cite{n1}  when $\F=\RR$ with the trivial flat
structure. Examples described in \refss{example}, show that using of
flat vector bundles $\F$ gives sharper estimates in general.

On the other hand, \refe{novin''} also generalizes the Morse type
inequalities obtained by S.P.Novikov in \cite{n3}, using Bloch
homology (which correspond to the case, when $[\omega]=0$ in
\refe{novin''}).

\ssec{example}{Examples} Here we describe examples, where the
Novikov numbers twisted by a flat vector bundle $\F$ (as defined
above) give greater values (and thus stronger inequalities) than the
usual Novikov numbers (where $\F=\RR \text{ or } {\CC}$,
cf. \cite{n1,p2}).

Let $k\subset S^3$ be a smooth knot and let the 3-manifold $X$ be the
result of $1/0$-surgery on $S^3$ along $k$. Note that the
one-dimensional homology group of $X$ is infinite cyclic and thus for
any complex number $\eta\in{\CC}$, $\eta\ne 0$, there is a complex
flat line bundle $\F_\eta$ over $X$ such that the monodromy with
respect to the generator of $H_1(X)$ is $\eta$.
By a choice of the knot $k$ and the number $\eta\in{\CC}^\ast$, we may
make the group $H^1(X,\F_\eta)$ arbitrarily large, while
$H^1(X,{\CC})$ is always one-dimensional.

Consider now the 3-manifold $M$ which is the connected sum
$
	M\ =\ X\ \#\ (S^1\times S^2).
$
Thus $M=X_+\cup X_-$,  where $X_+\cap X_-=S^2$,
$X_+=X-\text\{disk\}$ and $X_-=(S^1\times S^2)-\text\{disk\}$.
Suppose $\F$ is a flat complex line bundle $\F$ over $M$ such that its
restriction to $X_+$ is isomorphic to $\F_\eta|_{X_+}$. Consider the
class $\xi\in H^1(X,\RR)$ such that its restriction onto $X_+$ is
trivial and its restriction to $X_-$ is the generator.

By using the Mayer-Vietoris sequence, we show that
$\beta_1(\xi,\F)\ =\ \dim_{{\CC}} H^1(X,\F_\eta)$.
As we noticed above, this number can be arbitrarily large, while
$\dim_{\CC} H^1(M,{\CC})=2$.

\ssec{method}{Sketch of the proof of \reft{main}}
Our proof of \reft{main} is based on a slight modification of the Witten
deformation  \cite{wi}  suggested by Bismut \cite{bis}
in his proof of the degenerate Morse inequalities of Bott.
However  our proof is rather different from \cite{bis} even in the case
$[\w]=0$. We entirely avoid the probabilistic
analysis of the heat kernels, which is the most difficult part of \cite{bis}.
Instead, we give an explicit estimate on
the number of the "small" eigenvalues of the deformed Laplacian.
We now will explain briefly the main steps of the proof.

Let $U$ be a small tubular neighborhood of $C$ in $M$. We  identify $U$ with
a neighborhood of the zero section in the normal bundle $\nu(C)$. Fix
an affine connection on $\nu(C)$. This connection defines a bigrading
$$
	\Ome^\b(M,\F)=\bigoplus \Ome^{i,j}(M,\F),
$$
where $\Ome^{i,j}(M,\F)$ is the space of forms having degree $i$ in
the horizontal direction and degree $j$ in the vertical direction.
For $s\in \RR$, let $\tau_s$ be the map from $\Ome^\b(M,\F)$ to itself
which sends $\alp\in \Ome^{i,j}(M,\F)$ to $s^j\alp$.

Following Bismut, we introduce a 2-parameter deformation
\eq{n-talp}
	\n_{t,\alp}:\Ome^\b(M,\F)\to \Ome^{\b+1}(M,\F),
		\qquad t,\alp \in \RR
\end{equation}
of the covariant derivative $\n$, such that, for large values of
$t,\alp$ the Betti numbers of the deformed de Rham complex
$\big(\Ome^\b(M,\F),\n_{t,\alp}\big)$ are equal to the Novikov numbers
$\bet_p(\xi,\F)$.  Let $e(\w):\Ome^\b(M,\F)\to \Ome^{\b+1}(M,\F)$
denote the external multiplication by \w. Then on $U$ the deformation
\refe{n-talp} is given by
\eq{near}
	\n_{t,\alp}=(\t)^{-1}\circ\big(\n+t\alp e(\w)\big)\circ\t,
\end{equation}
while outside of some larger  neighborhood $V\supset U$, we have
\/ $\n_{t,\alp}=\n+t\alp e(\w)$.

There exists a unique function $f:U\to \RR$ such that $df=\w$ and
$f_{|_C}=0$. By the parameterized Morse lemma there exist an Euclidean
metric $h^{\nu(C)}$ on $\nu(C)$ such that $\nu(C)$ decomposes into an
orthogonal direct sum $\nu(C)=\nu^+(C)\oplus \nu^-(C)$ and
if $(y^+,y^-)\in U$, then $f(y)=\frac{|y^+|^2}{2}-\frac{|y^-|^2}{2}$.  Fix an
arbitrary Riemannian metric $g^C$ on $C$. The metrics $h^{\nu(C)}, g^C$
define naturally a Riemannian metric $g^{\nu(C)}$ on $\nu(C)$ (here we
consider $\nu(C)$ as a non-compact manifold).

Let $g^M$ be any Riemannian metric on $M$ whose restriction to $U$ is
equal to $g^{\nu(C)}$.  We also choose a Hermitian metric $h^\F$ on
$\F$.  Let us denote by $\Del_{t,\alp}$ the Laplacian associated with
the differential \refe{n-talp} and with the metrics $g^M,h^\F$.

Fix $\alp>0$ sufficiently large. It turns out that, when $t\to\infty$,
the eigenfunctions of $\Del_{t,\alp}$ corresponding to ``small''
eigenvalues localize near the critical points set $C$ of \w. Hence,
the number of the "small" eigenvalues of $\Del_{t,\alp}$ may be
calculated by means of the restriction of $\Del_{t,\alp}$ on $U$.  We
are led, thus, to study of a certain Laplacian on $\nu(C)$. The latter
Laplacian may be decomposed as $\bigoplus_Z\Del^Z_{t,\alp}$ where the
sum ranges over all connected components of $C$ and $\Del^Z_{t,\alp}$
is a Laplacian on the normal bundle $\nu(Z)=\nu(C)_{|_Z}$ to $Z$. We
denote by $\Del^{Z,p}_{t,\alp} \ (p=0,1,2,\dots)$ the restriction of
$\Del^Z_{t,\alp}$ on the space of $p$-forms.

It follows from  \refe{near}, that the spectrum of
$\Del^Z_{t,\alp}$ does not depend on $t$. Moreover, if $\alp>0$ is
sufficiently large, then
\eq{bismut'}
	\dim\Ker\Del^{Z,p}_{t,\alp}=
		\dim H^{p-\ind(Z)}(Z,\F_{|_Z}\otimes o(Z)).
\end{equation}
In the case when $\F$ is a trivial line bundle,
\refe{bismut'} is proven by Bismut \cite[Theorem 2.13]{bis}.

Let $E^\b_{t,\alp} \ (p=0,1\nek n)$ be the subspace of $\Ome^\b(M,\F)$
spanned by the eigenvectors of $\Del_{t,\alp}$ corresponding to the
``small'' eigenvalues. The cohomology of the deformed de Rham complex
$\big(\Ome^\b(M,\F),\n_{t,\alp}\big)$ may be calculated as the
cohomology of the subcomplex $\big(E^\b_{t,\alp},\n_{t,\alp}\big)$.

Using the method of \cite{sh2}, we show that, if the parameters $t$
and $\alp$ are large enough, then
\eq{k=h'}
	\dim E^p_{t,\alp}= \sum_{Z}\dim\Ker\Del^{Z,p}_{t,\alp},
\end{equation}
where the sum ranges over all connected components $Z$ of $C$. The
\reft{main} follows now from \refe{bismut'},\refe{k=h'} by standard
arguments (cf. \cite{bott2}).
\begin{Rem}In \cite{hs3}, Helffer and Sj\"ostrand  gave a
  very elegant analytic proof of the degenerate Morse inequalities of
  Bott. Though they also used the ideas of \cite{wi}, their method is
  completely different from \cite{bis}.  It is not clear if this method may
  be applied to the case $\xi\not=0$.
\end{Rem}

\ssec{L2}{$L^2$ generalization} Our \reft{main}, combined with the
results of W.L\"uck \cite{lu}, gives the following $L^2$ version of
the Novikov-Bott inequalities (5). Recall that $L^2$ generalization of
the usual Morse inequalities for Morse functions were obtained first
by S.P.Novikov and M.A.Shubin in \cite{ns}. $L^2$-version of Novikov
inequalities for 1-forms (allowing only isolated critial points) is
considered in a recent preprint \cite{MS} of V.Mathai and
M.Shubin. They use different technique and their assumptions do not
require resudual finiteness.

Let $\pi$ be a countable residually finite group and let $\calN(\pi)$ denote
the von Neumann algebra of $\pi$ acting on the Hilbert space $l^2(\pi)$
from the left and commuting with the standard action of $\pi$
on $l^2(\pi)$ from the right. The algebra $\calN(\pi)$ is supplied with the
canonical finite trace and all von Neumann dimensions later
will be understood with
respect to this trace.

Suppose that a flat bundle $\calL^{\pi}$ of Hilbert
$\calN(\pi)$-modules $l^2(\pi)$ over a closed manifold $M$ is given.
(Here $\pi$ is not necessarily the fundamental group of $M$).  Any
such bundle can be constructed by the standard construction from a
representation of the fundamental group of $M$ into $\pi$. Let $\F$
denote a finite dimensional flat vector bundle over $M$ as above. The
tensor product $\calL^{\pi}\otimes \F$ (the tensor product taken over
$\CC$) is again a bundle of Hilbert $\calN(\pi)$-modules over $M$.

Let \w\ be a  closed real valued 1-form  on $M$, which is
non-degenetrate in the sense of Bott. It determines a family of flat
bundles $\F_t$ as in \refss{novnum}.
Then  there exists a countable subset $S\in \RR$ (the set of
{\em jump points}) such that the von Neumann dimension
$$
	\dim_{\calN(\pi)}H^i_{(2)}(M,\calL^{\pi}\otimes\F_t)
$$
is constant for $t\notin S$ and it jumps up for $t\in S$. This fact
follows from Theorem 0.1 of L\"uck \cite{lu}. We will define {\it the von
Neumann - Novikov numbers} $\beta_i(\xi,\calL^\pi\otimes\F)$ as the
value of $\dim_{\calN(\pi)}H^i_{(2)}(M,\calL^{\pi}\otimes\F_t)$ for
$t\notin S$.  This value clearly depends only on the cohomology class
$\xi\in H^1(M,\RR)$ of $\omega$. Define {\it the von Neumann - Novikov
polynomial}
$$
	\calN_{\xi,\calL^{\pi}\otimes\F}(\lambda)\ = \
		\sum \lambda^i\beta_i(\xi,\calL^{\pi}\otimes\F).
$$
For any component $Z$ of the set of critical points $C$ of $\ome$ define
the following {\it von Neumann - Poincar\'e polynomial}
$$
    \calP_{Z, \calL^{\pi}\otimes \F}(\lam)\ =\  \sum
	\lam^i\dim_{\calN(\pi)}H^i_{(2)}(Z,\calL^{\pi}_{|_Z}
			\otimes\F_{|_Z}\otimes o(Z)),
$$
and then the {\it von Neumann - Morse counting polynomial}
$$
	\calM_{\omega,\calL^{\pi}\otimes\F}(\lambda)\ =
		\sum_{Z} \lam^{\ind(Z)}\calP_{Z, \calL^{\pi}\otimes \F}(\lam),
$$
the sum is taken over the set of connected components $Z$ of  $C$.

\th{L2} There exists a polynomial $\calQ(\lam)$ with real
non-negative coefficients, such that
$$
	\calM_{\ome,\calL^{\pi}\otimes\F}(\lam)-
		\calN_{\xi,\calL^{\pi}\otimes\F}(\lambda) =
					(1+\lambda)\calQ(\lambda).
$$
\eth

The proof is based on  theorem (0.1) of L\"uck \cite{lu} and \reft{main}.
Let's briefly indicate the main points.

Let $\pi\supset \Gamma_1\supset\Gamma_2\supset\dots $ be a sequence of
subgroups of finte index in $\pi$ having trivial intersection. The flat
$\calN(\pi)$ bundle $\calL^{\pi}$ is constructed by means of a representation
$\psi :\pi_1(M)\to \pi$. Let $\psi_m :\pi_1(M)\to \pi/\Gamma_m$ denote the
composition of $\psi$ with the reduction modulo $\Gamma_m$. The representation
$\psi_m$ determines
a flat vector bundle $\calL^{\pi}_m$ whose fiber is the group ring
$\CC[\pi/\Gamma_m]$ for any $m$. Slightly generalizing theorem (0.1) of
L\"uck \cite{lu}, we obtain that
$$
	\dim_{\calN(\pi)}H^i_{(2)}(M,\calL^{\pi}\otimes\F)=
	  \lim_{m\to\infty}
	    |\pi/\Gamma_m|^{-1}\dim_{\CC}H^i(M,\calL^{\pi}_m\otimes\F)
$$
This allows to approximate the von Neumann - Novikov polynomial
$\calN_{\xi,\calL^{\pi}\otimes\F}(\lam)$ by the polynomials
$$
	|\pi/\Gamma_m|^{-1}\calN_{\xi,\calL^{\pi}_m\otimes\F}(\lambda).
$$
Similarly, the von Neumann - Morse polynomial
$\calM_{\ome,\calL^{\pi}\otimes\F}(\lam)$ is approximated by the
polynomials
$|\pi/\Gam_m|^{-1}\calM_{\ome,\calL^{\pi}_m\otimes\F}(\lam)$. Application of
\reft{main} then finishes the proof.


\end{document}